\input amstex
\input amsppt.sty
\magnification=\magstep1
\hsize=30truecc
\vsize=22.2truecm
\baselineskip=16truept
\TagsOnRight
\nologo
\pageno=1
\topmatter
\def\N{\Bbb N}
\def\Z{\Bbb Z}

\def\l{\left}
\def\r{\right}
\def\b{\bigg}

\def\({\b(}
\def\[{\b[}
\def\){\b)}
\def\]{\b]}

\def\t{\text}
\def\f{\frac}
\def\mo{\roman{mod}}

\def\eq{\equiv}

\def\ls{\leqslant}

\def\al{\alpha}

\def\Proof{\noindent{\it Proof}}
\def\Remark{\medskip\noindent{\it Remark}}
\def\Ack{\noindent {\bf Acknowledgment}}
\hbox {J. Number Theory 129(2009), no.\,4, 964-969.}
\bigskip
\title Mixed sums of squares and triangular numbers (III)\endtitle
\author Byeong-Kweon Oh$^1$ and Zhi-Wei Sun$^2$\endauthor
\leftheadtext{Byeong-Kweon Oh and Zhi-Wei Sun}
\affil $^1$Department of Applied Mathematics, Sejong University
\\ Seoul, 143-747, Republic of Korea
\\bkoh\@sejong.ac.kr
\medskip
$^2$Department of Mathematics, Nanjing University\\
 Nanjing 210093, People's Republic of China\\zwsun\@nju.edu.cn
\\{\tt http://math.nju.edu.cn/$\sim$zwsun}
\endaffil
\thanks 2000 {\it Mathematics Subject Classification}.\,Primary 11E25;
Secondary 05A05, 11D85, 11P99, 11Y11.
\newline\indent
The first author
is supported by the Korea Research Foundation Grant (KRF-2008-314-C00004) funded by the Korean Government.
\newline\indent
The second author is responsible for communications, and supported by the National Natural Science Foundation
(grant 10871087) of People's Republic of China.
\endthanks
\abstract In this paper we confirm a conjecture of Sun which states that each positive integer
is a sum of a square, an odd square and a triangular number. Given any positive integer $m$,
we show that $p=2m+1$ is a prime congruent to 3 modulo 4 if and only if
$T_m=m(m+1)/2$ cannot be expressed as a sum of two odd squares and a triangular number,
i.e., $p^2=x^2+8(y^2+z^2)$ for no odd integers $x,y,z$.
We also show that a positive integer cannot be written as a sum of an odd square and two triangular numbers
if and only if it is of the form $2T_m\ (m>0)$ with  $2m+1$ having no prime divisor congruent to 3 modulo 4.
\endabstract
\endtopmatter

\document

\heading{1. Introduction}\endheading

The study of expressing natural numbers as sums of squares has long
history. Here are some well-known classical results in number
theory.

 (a) (Fermat-Euler theorem) Any prime $p\eq 1\ (\mo\ 4)$
 is a sum of two squares of integers.

 (b) (Gauss-Legendre theorem, cf. [G, pp.\,38--49] or [N, pp.\,17--23])
 $n\in\N=\{0,1,2,\ldots\}$ can be written as a sum of three squares
 of integers if and only if $n$ is not of the form $4^k(8l+7)$ with $k,l\in\N$.

 (c) (Lagrange's theorem) Every $n\in\N$ is a sum of four squares of integers.

Those integers $T_x=x(x+1)/2$ with $x\in\Z$ are called triangular numbers.
Note that $T_x=T_{-x-1}$ and $8T_x+1=(2x+1)^2$. In 1638 P. Fermat asserted
 that each $n\in\N$ can be written as a sum of three triangular numbers
 (equivalently, $8n+3$ is a sum of three squares of odd integers);
  this follows from the Gauss-Legendre theorem.

  Let $n\in\N$. As observed by L. Euler (cf. [D,
p.\,11]), the fact that $8n+1$ is a sum of three squares (of integers) implies
that $n$ can be expressed as a sum of two squares and a triangular
number. This is remarkable since there are infinitely many natural numbers
which cannot be written as a sum of three squares.
According to [D, p.\,24], E. Lionnet stated, and V. A.
Lebesgue [L] and M. S. R\'ealis [R] showed that $n$ is also a sum
of two triangular numbers and a square. In 2006 these two results were re-proved
by H. M. Farkas [F] via the theory of theta functions.
Further refinements of these results are summarized in the following theorem.

\proclaim{Theorem 1.0} {\rm (i) (B. W. Jones and G. Pall [JP])}
For every $n\in\N$, we can write $8n+1$ in the form $8x^2+32y^2+z^2$ with $x,y,z\in\Z$,
i.e., $n$ is a sum of a square, an even square and a triangular number.

{\rm (ii) (Z. W. Sun [S07])} Any natural number
is a sum of an even square and two triangular numbers.
If $n\in\N$ and $n\not=2T_m$ for any $m\in\N$, then $n$ is also a sum of an odd square and two triangular numbers.

{\rm (iii) (Z. W. Sun [S07])} A positive integer is a sum of an odd square, an even square and a triangular
number unless it is a triangular number $T_m\ (m>0)$ for which all prime divisors of $2m+1$ are congruent to $1$ mod $4$.
\endproclaim

We mention that Jones and Pall [JP] used the theory of ternary quadratic forms
and Sun [S07] employed some identities on $q$-series. Motivated by Theorem 1.0(iii) and the fact that
every prime $p\eq1\ (\mo\ 4)$ is a sum of an odd square and an even square, the second author [S09]
conjectured that each natural number $n\not=216$ can be written in the form
$p+T_x$ with $x\in\Z$, where $p$ is a prime or zero.
Sun [S09] also made a general conjecture which states that for any $a,b\in\N$ and $r=1,3,5,\ldots$ all sufficiently large
integers can be written in the form $2^ap+T_x$ with $x\in\Z$, where $p$ is either zero or a prime
congruent to $r$ mod $2^b$.

In [S07] Sun investigated what kind of mixed sums $ax^2+by^2+cT_z$
or $ax^2+bT_y+cT_z$ (with $a,b,c\in\Z^+=\{1,2,3,\ldots\}$) represent all natural numbers, and left two conjectures
in this direction. In [GPS] S. Guo, H. Pan and Sun proved Conjecture 2 of [S07]. Conjecture 1 of Sun [S07]
states that any positive integer $n$ is a sum of a square, an odd square and a triangular number, i.e.,
$n-1=x^2+8T_y+T_z$ for some $x,y,z\in\Z$.

 In this paper we prove Conjecture 1 of Sun [S07] and some other results concerning mixed sums of squares and triangular numbers.
 Our main result is as follows.

\proclaim{Theorem 1.1} {\rm (i)} Each positive integer is a sum of a square, an odd square and a triangular number.
A triangular number $T_m$ with $m\in\Z^+$ is a sum of two odd squares and a triangular number
if and only if $2m+1$ is not a prime congruent to $3$ mod $4$.

{\rm (ii)} A positive integer cannot be written as a sum of an odd square and two triangular numbers
if and only if it is of the form $2T_m\ (m\in\Z^+)$ with $2m+1$ having no prime divisor congruent to $3$ mod $4$.
\endproclaim

\Remark\ 1.1. In [S09] the second author conjectured that if a positive integer is not a triangular number
then it can be written as a sum of two odd squares and a triangular number unless it is among the following
$25$ exceptions:
$$\align &4,\ 7,\ 9,\ 14,\ 22,\ 42,\ 43,\ 48,\ 52,\ 67,\ 69,\ 72,\ 87,\ 114,
\\&144,\ 157,\ 159,\ 169,\ 357,\ 402,\ 489,\ 507,\ 939,\ 952,\ 1029.
\endalign$$

\smallskip

Here is a consequence of Theorem 1.1.

\proclaim{Corollary 1.1} {\rm (i)} An odd integer $p>1$ is a prime congruent to $3$ mod $4$ if and only if
$p^2=x^2+8(y^2+z^2)$ for no odd integers $x,y,z$.

{\rm (ii)} Let $n>1$ be an odd integer. Then all prime divisors of $n$ are congruent to $1$ mod $4$,
if and only if $n^2=x^2+4(y^2+z^2)$ for no odd integers $x,y,z$.
\endproclaim

\Remark\ 1.2. In number theory there are very few simple characterizations of primes such as Wilson's theorem.
Corollary 1.1(i) provides a surprising new criterion for primes congruent to 3 mod 4.

\medskip
In the next section we will prove an auxiliary theorem.
Section 3 is devoted to our proofs of Theorem 1.1 and Corollary 1.1.

\heading{2. An auxiliary theorem}\endheading

 In this section we prove the following auxiliary result.

\proclaim{Theorem 2.1} Let $m$ be a positive integer.

{\rm (i)} Assume that $p=2m+1$ be a prime congruent to $3$ mod $4$. Then
$T_m$ cannot be written in the form $x^2+y^2+T_z$ with $x,y,z\in\Z$,
$x^2+y^2>0$ and $x\eq y\ (\mo\ 2)$. Also, $2T_m$ is not a sum of a positive even square and
two triangular numbers.

{\rm (ii)} Suppose that
all prime divisors of $2m+1$ are congruent to $1$ mod $4$.
Then
$T_m$ cannot be written as  a sum of an odd square, an even square and a triangular number.
Also, $2T_m$ is not a sum of an odd square and two triangular numbers.
\endproclaim

To prove Theorem 2.1 we need the following result due to Hurwitz.

\proclaim{Lemma 2.1 {\rm (cf. [D, p.\,271] or [S07, Lemma 3])}} Let $n$ be a positive odd integer,
and let $p_1,\ldots,p_r$ be all the distinct prime divisors of $n$
congruent to $3$ mod $4$. Write $n=n_0\prod_{0<i\ls r}p_i^{\alpha_i}$, where $n_0,\al_1,\ldots,\al_r\in\Z^+$
and $n_0$ has no prime divisors congruent to $3$ mod $4$. Then
$$|\{(x,y,z)\in\Z^3:\,x^2+y^2+z^2=n^2\}|=6n_0\prod_{0<i\ls r}\l(p_i^{\al_i}+2\f{p_i^{\al_i}-1}{p_i-1}\r).$$
\endproclaim

 As in [S07], for $n\in\N$ we define
  $$r_0(n)=|\{(x,y,z)\in\Z\times\N\times\N:\, x^2+T_y+T_z=n\
 \t{and}\ 2\mid x\}|$$
 and
$$r_1(n)=|\{(x,y,z)\in\Z\times\N\times\N:\, x^2+T_y+T_z=n\
 \t{and}\ 2\nmid x\}|.$$

 By p.\,108 and Lemma 2 of Sun [S07], we have the following lemma.

 \proclaim{Lemma 2.2 {\rm ([S07])}} For $n\in\N$ we have
 $$|\{(x,y,z)\in\Z\times\Z\times\N:\, x^2+y^2+T_z=n\
 \t{and}\ x\eq y\ (\mo\ 2)\}|=r_0(2n)$$
 and
$$|\{(x,y,z)\in\Z\times\Z\times\N:\, x^2+y^2+T_z=n\
 \t{and}\ x\not\eq y\ (\mo\ 2)\}|=r_1(2n).$$
 Also, $r_0(2T_m)-r_1(2T_m)=(-1)^m(2m+1)$ for every $m\in\N$.
 \endproclaim

\medskip

\noindent{\it Proof of Theorem 2.1}.
By Lemma 2.2,
$$\align &r_0(2T_m)+r_1(2T_m)
\\=&|\{(x,y,z)\in\Z\times\Z\times\N:\, x^2+y^2+T_z=T_m\}|
\\=&\f12|\{(x,y,z)\in\Z^3:\, 8x^2+8y^2+(8T_z+1)=8T_m+1\}|
\\=&\f12|\{(x,y,z)\in\Z^3:\, 4(x+y)^2+4(x-y)^2+(2z+1)^2=(2m+1)^2\}|
\\=&\f12|\{(u,v,z)\in\Z^3:\, 4(u^2+v^2)+(2z+1)^2=(2m+1)^2\}|
\\=&\f16|\{(x,y,z)\in\Z^3:\, x^2+y^2+z^2=(2m+1)^2\}|.
\endalign$$

(i) As $p=2m+1$ is a prime congruent to 3 mod 4,
by Lemma 2.1 and the above we have
$$r_0(2T_m)+r_1(2T_m)=p+2.$$
On the other hand,
$$r_0(2T_m)-r_1(2T_m)=(-1)^{m}(2m+1)=-p$$
by Lemma 2.2.
So
$$2r_0(2T_m)=r_0(2T_m)+r_1(2T_m)+(r_0(2T_m)-r_1(2T_m))=p+2-p=2.$$
Therefore
$$|\{(x,y,z)\in\Z\times\Z\times\N:\, x^2+y^2+T_z=T_m\
 \t{and}\ 2\mid x-y\}|=r_0(2T_m)=1$$
 and also
 $$|\{(x,y,z)\in\Z\times\N\times\N:\, x^2+T_y+T_z=2T_m\
 \t{and}\ 2\mid x\}|=r_0(2T_m)=1.$$
 Since $T_m=0^2+0^2+T_m$ and $2T_m=0^2+T_m+T_m$, the desired results follow immediately.

 \medskip

 (ii) As all prime divisors of $2m+1$ are congruent to 1 mod 4,
 we have
$$r_0(2T_m)+r_1(2T_m)=\f16|\{(x,y,z)\in\Z^3:\, x^2+y^2+z^2=(2m+1)^2\}|=2m+1$$
in view of Lemma 2.1. Note that $m$ is even since $2m+1\eq1\ (\mo\ 4)$.
By Lemma 2.2,  $r_0(2T_m)-r_1(2T_m)=(-1)^{m}(2m+1)=2m+1$. Therefore
$$|\{(x,y,z)\in\Z\times\Z\times\N:\, x^2+y^2+T_z=T_m\
 \t{and}\ 2\nmid x-y\}|=r_1(2T_m)=0.$$
This proves part (ii) of Theorem 2.1. \qed

\heading{3. Proofs of Theorem 1.1 and Corollary 1.1}\endheading

\proclaim{Lemma 3.1} Let $m\in\N$ with $2m+1=k(w^2+x^2+y^2+z^2)$ where $k,w,x,y,z\in\Z$.
Then
$$2T_m=k^2(wy+xz)^2+k^2(wz-xy)^2+2T_v\quad\t{for some}\ v\in\Z.$$
\endproclaim
\Proof. Write the odd integer $k(w^2+x^2-(y^2+z^2))$ in the form $2v+1$. Then
$$\align 8T_m+1=&(2m+1)^2=k^2(w^2+x^2+y^2+z^2)^2
\\=&(2v+1)^2+4k^2(w^2+x^2)(y^2+z^2)
\\=&8T_v+1+4k^2((wy+xz)^2+(wz-xy)^2)
\endalign$$
and hence
$$2T_m=2T_v+k^2(wy+xz)^2+k^2(wz-xy)^2.$$
This concludes the proof. \qed

\medskip
\noindent {\it Proof of Theorem 1.1}. (i) In view of Theorem
1.0(iii), it suffices to show the second assertion in part (i).

Let $m$ be any positive integer. By Theorem 2.1(i), if $2m+1$ is a prime congruent to 3 mod 4
then $T_m$ cannot be written as a sum of two odd squares and a triangular number.

Now assume that $2m+1$ is not a prime congruent to 3 mod 4.
Since the product of two integers congruent to 3 mod 4
is congruent to 1 mod 4, we can write $2m+1$ in the form $k(4n+1)$
with $k,n\in\Z^+$.

Set $w=1+(-1)^n$. Observe that $4n+1-w^2$ is a positive integer congruent to 5 mod
8. By the Gauss-Legendre theorem on sums of three squares, there are
integers $x,y,z$ with $x$ odd such that $4n+1-w^2=x^2+y^2+z^2$.
Clearly both $y$ and $z$ are even. As $y^2+z^2\eq 4\ (\mo\ 8)$, we
have $y_0\not\eq z_0\ (\mo\ 2)$ where $y_0=y/2$ and $z_0=z/2$.

Since $2m+1=k(w^2+x^2+y^2+z^2)$, by Lemma 3.1 there is an integer $v$ such that
$$2T_{m}=2T_v+k^2(wy+xz)^2+k^2(wz-xy)^2.$$
Thus
$$\align T_{m}=&T_v+2(kwy_0+kxz_0)^2+2(kwz_0-kxy_0)^2
\\=&T_v+(kwy_0+kxz_0+(kwz_0-kxy_0))^2
\\&+(kwy_0+kxz_0-(kwz_0-kxy_0))^2.
\endalign$$
As $w$ is even and $kxy_0\eq y_0\not\eq z_0\eq kxz_0\ (\mo\ 2)$,
we have $$kwy_0+kxz_0\pm(kwz_0-kxy_0)\eq1\ (\mo\ 2).$$
Therefore $T_{m}-T_v$ is a sum of two odd squares.

\medskip

(ii) In view of Theorem 1.0(ii) and Theorem 2.1(ii), it suffices to show that
if $2m+1\ (m\in\Z^+)$ has a prime divisor congruent to 3 mod 4 then
$2T_m$ is a sum of an odd square and two triangular numbers.

Suppose that $2m+1=k(4n-1)$ with $k,n\in\Z^+$. Write $w=1+(-1)^n$.
Then $4n-1-w^2$ is a positive integer congruent to 3 mod 8.
By the Gauss-Legendre theorem on sums of three squares, there are integers
$x,y,z$ such that $4n-1-w^2=x^2+y^2+z^2$. Clearly $x\eq y\eq z\eq1\ (\mo\ 2)$
and $2m+1=k(w^2+x^2+y^2+z^2)$.
By Lemma 3.1, for some $v\in\Z$ we have
$$2T_{m}=k^2(wy+xz)^2+k^2(wz-xy)^2+2T_v.$$
Let $u=kwz-kxy$. Then
$$T_{v+u}+T_{v-u}=\f{(v+u)^2+(v-u)^2+(v+u)+(v-u)}2=u^2+2T_v.$$
Thus
$$2T_{m}=(kwy+kxz)^2+T_{v+u}+T_{v-u}.$$
Note that $kwy+kxz$ is odd since $w$ is even and $k,x,z$ are odd.

\smallskip
Combining the above we have completed the proof of Theorem 1.1. \qed

\medskip
\noindent{\it Proof of Corollary 1.1}. (i) Let $m=(p-1)/2$. Observe that
$$\align &T_m=T_x+(2y+1)^2+(2z+1)^2
\\\iff &p^2=8T_m+1=(2x+1)^2+8(2y+1)^2+8(2z+1)^2.
\endalign$$
So the desired result follows from Theorem 1.1(i).

(ii) Let $m=(n-1)/2$. Clearly
$$\align &2T_m=T_x+T_y+(2z+1)^2
\\\iff&2n^2=16T_m+2=(2x+1)^2+(2y+1)^2+8(2z+1)^2
\\\iff&n^2=(x+y+1)^2+(x-y)^2+4(2z+1)^2.
\endalign$$
So $2T_m$ is a sum of an odd square and two triangular numbers if
and only if $n^2=x^2+(2y)^2+4z^2$ for some odd integers $x,y,z$.
(If $x$ and $z$ are odd but $y$ is even, then $x^2+(2y)^2+4z^2\eq
5\not\eq n^2\ (\mo\ 8)$.) Combining this with Theorem 1.1(ii) we
obtain the desired result. \qed

\medskip

\Remark\ 3.1. We can deduce Corollary 1.1 in another way by using some known results (cf. [E], [EHH] and [SP])
in the theory of ternary quadratic forms, but this approach involves many sophisticated concepts.
\medskip

\Ack. The authors are grateful to the referee for his/her helpful comments.


\bigskip

\widestnumber\key{EHH}
\Refs


\ref\key D\by L. E. Dickson\book
History of the Theory of Numbers, {\rm Vol. II}
\publ AMS Chelsea Publ., 1999\endref

\ref\key E\by A. G. Earnest\paper Representation of spinor exceptional integers by ternary quadratic forms
\jour Nagoya Math. J.\vol 93\yr 1984\pages 27--38\endref

\ref\key EHH\by A. G. Earnest, J. S. Hsia and D. C. Hung\paper
Primitive representations by spinor genera of ternary quadratic forms
\jour J. London Math. Soc.\vol 50\yr 1994\pages 222--230\endref

\ref\key F\by H. M. Farkas\paper Sums of squares and triangular numbers
\jour Online J. Anal. Combin.\vol 1\yr 2006\pages \#1, 11 pp. (electronic)\endref

\ref\key G\by E. Grosswald\book Representation of Integers as Sums of Squares
\publ Springer, New York, 1985\endref

\ref\key GPS\by S. Guo, H. Pan and Z. W. Sun\paper
Mixed sums of squares and triangular numbers (II)
\jour Integers\vol 7\yr 2007\pages \#A56, 5pp (electronic)\endref

\ref\key JP\by B. W. Jones and G. Pall
\paper Regular and semi-regular positive ternary quadratic forms
\jour Acta Math.\vol 70\yr 1939\pages 165--191\endref

\ref\key L\by V. A. Lebesque\paper Questions 1059,1060,1061 (Lionnet)
\jour Nouv. Ann. Math.\vol 11\yr 1872\pages 516--519\endref

\ref\key N\by M. B. Nathanson\paper Additive Number Theory: The
Classical Bases \publ Grad. Texts in Math., vol. 164, Springer,
New York, 1996\endref

\ref\key R\by M. S. R\'ealis\paper Scolies pour un th\'eoreme d'arithm\'etique
\jour Nouv. Ann. Math.\vol 12\yr 1873\pages 212--217\endref

\ref\key SP\by R. Schulze-Pillot\paper Darstellung durch Spinorgeschlechter ternarer quadratischer Formen {\rm (German)}
\jour J. Number Theory\vol 12\yr 1980\pages 529--540\endref

\ref\key S07\by Z. W. Sun\paper Mixed sums of
squares and triangular numbers \jour Acta Arith. \vol 127\yr 2007\pages 103--113\endref

\ref\key S09\by Z. W. Sun\paper On sums of primes and triangular
numbers \jour Journal of Combinatorics and Number Theory \vol 1\yr
2009\pages 65--76. {\tt http://arxiv.org/abs/0803.3737}\endref

\endRefs
\enddocument